\documentclass{article}

\usepackage{PRIMEarxiv}

\usepackage[utf8]{inputenc} 
\usepackage[T1]{fontenc}    
\usepackage{hyperref}       
\usepackage{url}            
\usepackage{booktabs}       
\usepackage{amsfonts}       
\usepackage{nicefrac}       
\usepackage{microtype}      
\usepackage{lipsum}
\usepackage{fancyhdr}       
\usepackage{graphicx}       

\usepackage{adjustbox}
\usepackage{algorithmic}
\usepackage{algorithm}
\usepackage{amsmath,amsfonts,amssymb,array}
\usepackage{bbm}
\usepackage{booktabs}
\usepackage{cancel}
\usepackage{caption}
\usepackage{cite}
\usepackage{color,colortbl}
\usepackage{epsfig}
\usepackage{epstopdf}
\usepackage{float}
\usepackage{graphicx}
\usepackage{lipsum}
\usepackage{lineno}
\usepackage{longtable}
\usepackage{mathtools}
\usepackage{multicol}
\usepackage{multirow}
\usepackage{pdflscape}
\usepackage{relsize}
\usepackage{rotating}
\usepackage{soul}
\usepackage{subfig}
\usepackage[dvipsnames,table]{xcolor}
\usepackage{tabularx}
\usepackage{tabulary}
\usepackage{textcomp}
\usepackage{tikz}
\usetikzlibrary{3d}
\usetikzlibrary{arrows,shapes,automata,backgrounds,petri,positioning}
\usetikzlibrary{decorations.pathmorphing}
\usetikzlibrary{decorations.shapes}
\usetikzlibrary{decorations.text}
\usetikzlibrary{decorations.fractals}
\usetikzlibrary{decorations.footprints}
\usetikzlibrary{shadows}
\usetikzlibrary{calc}
\usetikzlibrary{spy}

\emergencystretch=\maxdimen
\newcommand{\Hmax}{H^{\max}}

\newcommand{\Omax}{O^{\max}}
\newcommand{\Plos}{\mathbb{P}^\textrm{LoS}}
\newcommand{\Pnlos}{\mathbb{P}^\textrm{NLoS}}
\newcommand{\VCN}{V^\mathrm{CN}}
\newcommand{\VD}{V^\mathrm{D}}
\newcommand{\VC}{V^\mathrm{C}}
\newcommand{\VCS}{V^\mathrm{CS}}
\newcommand{\VF}{V^\mathrm{F}}
\DeclarePairedDelimiter\norm{\lVert}{\rVert}

\title{Communication-aware Drone Delivery Problem
}

\author{
  Cihan~Tugrul~Cicek \\
  Department of Industrial Engineering \\
  Atilim University \\
  Ankara, Turkey \\
  \texttt{cihan.cicek@atilim.edu.tr} \\
   \And
  Çağrı~Koç \\
  Department of Business Administration  \\
  Social Sciences University of Ankara \\
  Ankara, Turkey \\
  \texttt{cagri.koc@asbu.edu.tr} \\
  \And
  Hakan~Gultekin \\
  Department of Mechanical and Industrial Engineering  \\
  Sultan Qaboos University \\
  Muscat, Oman \\
  \texttt{hgultekin@squ.edu.om} \\
  \And
  G\"une\c s~Erdo\u gan \\
  School of Management \\
  University of Bath \\
  Bath, UK \\
  \texttt{ge277@bath.ac.uk} \\
}

\begin{document}
\maketitle

\begin{abstract}
The drone delivery problem (DDP) has been introduced to include aerial vehicles in last-mile delivery operations to increase efficiency. However, the existing studies have not incorporated the communication quality requirements of such a delivery operation. This study introduces the Communication-aware DDP (C-DDP), which incorporates handover and outage constraints. In particular, any trip of a drone to deliver a customer package must require less than a certain number of handover operations and cannot exceed a predefined outage duration threshold. The authors develop a Mixed Integer Programming (MIP) model to minimize the total flight distance while satisfying communication constraints as well as the time windows of customers. We present a Genetic Algorithm (GA) that can solve large instances, and compare its performance with an off-the-shelf MIP solver. Computational results show that the GA can outperform the MIP solver for solving larger instances and is a better option.
\end{abstract}

\keywords{Drone delivery \and genetic algorithm \and outage \and handover}

\section{Introduction}\label{sec:introduction}

In the past decades, rapid technological advances have resulted in the increasing use of \emph{Unmanned Aerial Vehicles} (UAVs) in diverse sectors. UAVs are utilized for package delivery, monitoring, and providing wireless communication coverage in three main application areas; agriculture, environment, and defense \cite{swain2007suitability, macrina2020drone}. Major companies such as Amazon, Google, and DHL have started to use UAVs in their parcel delivery process \cite{yoo2018drone}.

Many UAV routing problems have been studied in the literature. The stochastic UAV routing problem is addressed in \cite{chow2016dynamic}, which monitors and directs UAV traffic under uncertainty. The authors developed a policy by adapting the Bellman equation and using an approximate algorithm. The multi-trip drone routing problem is studied in \cite{cheng2020drone}, where the energy consumption of the drones is modeled as a nonlinear function of payload and distance traveled. The authors developed logical and subgradient cuts in the solution process to implement the complex convex energy function. \cite{yuan2021enhanced} developed a genetic algorithm to solve the heterogeneous UAV routing problem. The problem of scheduling and sequencing drone routes for medical item delivery is studied in \cite{ghelichi2021logistics}, where locations for charging stations as well as assignment of clinics to providers are determined. The \emph{Vehicle Routing Problem} (VRP) with drones that incorporates time windows for the customers and requires coordination between trucks and drones is addressed in \cite{kuo2021vehicle}. The authors provided a mathematical model to minimize the total travel time of trucks and proposed a Variable Neighborhood Search algorithm. The drone routing problem with recharging stops is studied in \cite{ermaugan2022learning}, which has many applications, including precision agriculture, search-and-rescue, and military surveillance. The authors developed a learning-based heuristic algorithm to solve the problem. We refer the interested readers to the survey papers \cite{otto2018optimization, khoufi2019survey, thibbotuwawa2020unmanned, rojas2021unmanned, poikonen2021future} for an in-detail exposition to the drone routing literature.

The UAV routing problem has also been studied in the context of wireless communications. \cite{zhang2018} addressed a cellular-enabled UAV communication system consisting of one UAV and multiple \emph{Ground Base Stations} (GBSs) to minimize mission completion time. \cite{zhang2019} studied the problem of minimizing the UAV mission completion time by optimizing its trajectory. The authors used a quality-of-connectivity constraint of the GBS-UAV link specified by a minimum received signal-to-noise ratio target. \cite{xie2021connectivity} studied the three-dimensional path planning problem for cellular-connected UAVs, incorporating the impact of 3D antenna radiation patterns. \cite{fakhreddine2019handover} discussed cell selection and handover measurements for an aerial drone in a suburban environment. \cite{cherif2021disconnectivity} studied the disconnectivity-aware and energy-efficient UAV trajectory planning with minimum handovers. \cite{bulut2018trajectory} studied the trajectory optimization for cellular-connected UAVs with a disconnection duration constraint.

The drone routing literature covers many different topics regarding parcel delivery \cite{repoussis2009arc, yan2019graph, rojas2021unmanned}. To the best of our knowledge, communication constraints have not been incorporated into drone routing problems. However, the communication performance between the drones and the core network is required to (i) avoid collisions, (ii) schedule recharging and repair operations, and (iii) improve operational efficiency in fleet management. We aim to fill this gap between the UAV-based package delivery and UAV communication requirements by introducing the \emph{Communication-aware Drone Delivery Problem} (C-DDP), which addresses how to determine optimal drone trips to deliver customer packages subject to communication quality constraints.

\begin{figure}[!htb]
    \centering
    \includegraphics[width=.5\linewidth]{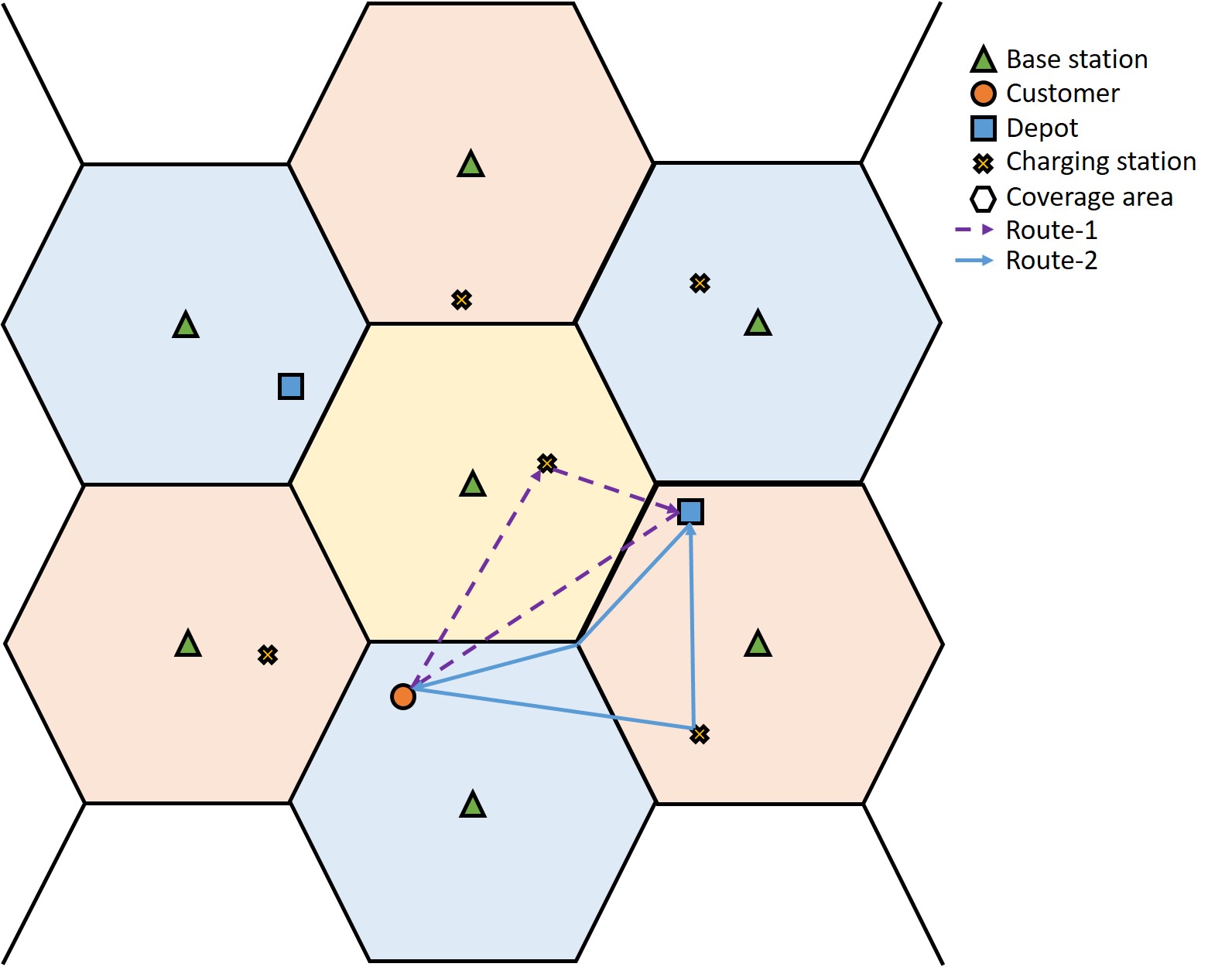}
    \caption{\it \textbf{Illustration of the C-DDP.} The dashed curves show a trip in which communication performance is disregarded. Consequently, the trip consists of straight flights between nodes to minimize the total flight distance. The solid curves show another trip in which the communication performance, e.g., number of handovers, is considered.}
    \label{fig:system_model}
\end{figure}

Figure~\ref{fig:system_model} presents an illustration of the C-DDP. Two alternative trips for a package delivery are depicted. In the first route, a drone is dispatched from the depot and follows a straight-line path to the customer, delivers the order, and visits a charging station before returning to the depot. This trip requires the drone to pass over three different coverage areas, risking handover. In the second route, the drone modifies its trajectory and visits only two coverage areas. Clearly, there is a trade-off between these two routes regarding flight distance and communication quality. The first route requires less flight distance but more handover operations, whereas the second route requires fewer handover operations, but more flight distance.

The rest of this paper is structured as follows. We present the problem statement, a model of the delivery system including the communication network and delivery network, as well as an illustrative example and a MIP formulation in Section \ref{sec:system_model}. We provide our solution method and computational results in Section \ref{sec:solution_approach}, followed by our conclusions in Section \ref{sec:conclusions}.

\section{System Model}\label{sec:system_model}

The C-DDP is defined on a complete directed graph $G=(V,A)$ with a set of vertices $V$ and a set of arcs $A$. The vertices consist of the communication nodes (CN), i.e., base stations, $\VCN = \{CN_0, CN_1, \ldots, CN_{n_{CN}-1}\}$, depots, $\VD = \{D_0, D_1, \ldots, D_{{n_D}-1}\}$, charging stations, $\VCS = \{CS_0, CS_1, \ldots, CS_{n_{CS}-1}\}$, and customers $\VC = \{C_0,C_1,\ldots,C_{{n_C}-1}\}$. As a result, we have $V=\VCN \cup \VD \cup \VCS \cup \VC$. The arcs consist of the pairs of vertices between which a drone can fly, i.e., $A = \{(i,j): i,j \in \VF\}$, where $\VF = \VD \cup \VCS \cup \VC$ denotes the flyable nodes.

\subsection{Communication Network}\label{subsec:communication_network}
The communication network consists of multiple fixed, ground CNs that have access to the core communication network. The entire service area is divided into multiple coverage areas corresponding to different CNs. During a trip, a drone may not be within the \emph{Line-of-Sight} (LoS) of a CN all the time; thus, it must establish \emph{reliable} communication links with the CNs to inform its status to a control center. The reliability in this context can be defined as experiencing a certain \emph{Quality-of-Service} (QoS) level, e.g., spectral efficiency (SE), higher than a predefined threshold during a trip.

Since a CN on the ground (i.e., GBS) is more capable of supporting technologies, we assume that \emph{Power Domain Non-Orthogonal Multiple Access} (PD-NOMA) and \emph{Orthogonal Frequency-Division Multiple Access} (OFDMA) are applied for the CNs and drones, respectively. Such technologies have already been proposed for aerial communications \cite{azizi2019,cicek2021}. We use the probabilistic \emph{pathloss} model of \cite{alhourani2014} that divides the users into two propagation groups. The first group is assumed to have a LoS connection with a probability of $\Plos$ with relatively low pathloss, while the second group can maintain a poor connection with a relatively high pathloss. We assume that the drones operate at a fixed altitude $H^\textrm{D}$ and establish a link with the closest CN in its corresponding coverage area.

Assume that a drone flies in the coverage area of a particular $CN_i \in \VCN$. Given the location of this CN ($g_i^\mathrm{CN}$) and the projected location of the drone on the ground ($g^\mathrm{D}$), the probability of LoS connection can be defined as
\begin{flalign}\label{eqn:probloss}
    \Plos_i & = \frac{1}{1+\alpha_1 e^{ -\alpha_2 \left( \frac{180}{\pi} \tan^{-1} \left( \frac{H^\textrm{D}}{\norm*{g_i^\mathrm{CN} - g^\mathrm{D}}} \right) - \alpha_1 \right)}},
\end{flalign}

\noindent where $\alpha_1$ and $\alpha_2$ are environment-dependent parameters. Accordingly, the \emph{non-line-of-sight} (NLoS) probability, i.e., probability of maintaining a poor connection, is defined as $\Pnlos_i = 1 - \Plos_i$. Then, the mean pathloss between the drone and $CN_i$, ($L_i$), can be formulated as,

{\small
\begin{flalign}\label{eqn:pathloss}
    \nonumber L_i & = 10 \alpha_3 \log_{10} \left( \frac{4\pi f_c }{c} \sqrt{\norm*{g_i^\mathrm{CN} - g^\mathrm{D}}^2 + (H^\textrm{d})^2} \right) \\
    & \qquad + \left( \mu_\mathrm{LoS} \Plos_i + \mu_\mathrm{NLoS} \Pnlos_i \right),
\end{flalign}
}

\noindent where $f_c$ is the carrier frequency in Hz, $c$ is the speed of light in m/s, $\alpha_3$ is the pathloss exponent, and $\mu_\mathrm{LoS}$ and $\mu_\mathrm{NLoS}$ are the environment-dependent average additional losses to the free-space propagation for LoS and NLoS connections, respectively \cite{zeng2019}.  Then, the signal-to-interference-plus-noise-ratio (SINR) at the drone receiving from $CN_i$, ($\rho_i$), is given as,
\begin{flalign}\label{eqn:sinr}
    \rho_i & = \dfrac{P_i 10^{-L_i \slash 10}}{\sigma^2 + \sum_{j \in \VCN,j \neq i}P_j 10^{-L_j \slash 10}}
\end{flalign}

\noindent where $P_i$ is the transmit power of $CN_i$ and $\sigma^2$ is the noise power at the drone receiver. Consequently, the SE when receiving from $CN_i$, ($\gamma_i$), is given as 
\begin{flalign}\label{eqn:spectral_efficiency}
    \gamma_i = \log_2(1+\rho_i).
\end{flalign}

Based on these definitions, we consider two different QoS measures. The first measure is the number of handovers. In particular, every time a flying drone establishes a new link to a new CN, a handover operation is required, i.e., transferring the bandwidth from one CN to the other CN.  We expect to have as few handover operations as possible during a trip. The second measure is the expected outage duration. The outage in this context can be defined as experiencing a SE that is less than a predefined SE threshold ($\Bar{\gamma}$).

More precisely, let $v_1,v_2 \in \VF$ denote two different nodes in the network. The number of handovers and outage probability for the straight path between these two nodes can be computed as the integrals on the path. However, discretization of the path, i.e., dividing into multiple line segments, has been shown as an effective approach for computing such measures in cellular connected trajectory design \cite{zhang2019,chen2020}. Thus, we adopt the same approach in this study. We divide the straight path between these two nodes into $R$ equal line segments, yielding $R+1$ points on this particular path (including $v_1$ and $v_2$). Let $b_r$ and $\gamma_r$ denote the closest CN and the SINR at point $r$ receiving from this closest CN, respectively. Then, the number of handovers, ($h_{{v_1}{v_2}}$), and the probability of outage on arc $(v_1,v_2)$, ($\mathbb{P}^\mathrm{O}_{{v_1}{v_2}}$), can be approximated as
\begin{flalign}
    h_{{v_1}{v_2}} & = \sum_{r=2}^{R+1} \mathbbm{1}(b_r \neq b_{r-1}), \label{eqn:number_of_handover}\\  
    \mathbb{P}^\mathrm{O}_{{v_1}{v_2}} & = \frac{1}{R+1} \sum_{r=1}^{R+1} \mathbbm{1}(\gamma_r < \Bar{\gamma}), \label{eqn:outage_probability}
\end{flalign}

\noindent where $\mathbbm{1}(z)$ is the indicator function that returns 1 if $z$ is \emph{True}, and 0 otherwise. 

Let $t_{{v_1}{v_2}}$ be the travel time between nodes $v_1$ and $v_2$. Then, the expected outage duration, ($o_{{v_1}{v_2}}$), on the arc $(v_1,v_2)$ can be computed as
\begin{flalign}\label{eqn:outage_duration}
    o_{{v_1}{v_2}} & = \mathbb{P}^\mathrm{O}_{{v_1}{v_2}} t_{{v_1}{v_2}}.
\end{flalign}

\subsection{Delivery Network}\label{subsec:delivery_network}
We consider a delivery network in which drones, operated by a control center, deliver the orders directly from depots to customers. We assume that the customer locations are known, the packages are not heavier than the payload capacity of drones, and a drone can deliver only one order on a trip. A trip is defined as a sequence of arcs starting from a depot node, visiting at most one customer node, and ending at a depot node.

The drones operate within a finite working day, and each drone is associated with a starting and an ending depot, which are not necessarily the same. A drone starts the day at its starting depot, can operate multiple trips, and must return to its ending depot. We assume that any depot can fulfill orders; thus, a drone is not enforced to return to the same depot from which it was loaded. To prevent crashing due to a depleted battery, a drone may visit a charging station on a trip to swap its battery with a fully charged one. We assume that the charging stations have infinite capacities, i.e., swapping starts immediately after a drone arrives at a charging station and requires a fixed and identical time at every station independent of the battery level of the arriving drone. Hence, the following constraints must be satisfied:
\begin{itemize}
    \item The first trip of a drone must start from its starting depot.
    \item The last trip of a drone must end at its ending depot.
    \item Each trip starts from a depot, ends at a depot, and can include at most one customer.
    \item Each customer is served exactly once by one drone.
    \item Each customer must be served within its time window.
\end{itemize}

\subsection{Problem Statement}\label{subsec:problem_statement}
We consider a fleet of homogeneous drones $U=\{U_1,U_2,\ldots,U_{n_U}\}$. Each customer $i \in \VC$ has a time window $[a_i,b_i]$ within which the delivery must be completed. Each flyable node $i \in \VF$ is associated with an operation time $w_i$ denoting the time that a drone must spend. This includes loading a package or battery swapping at a depot node, unloading a package at a customer node, or swapping at a charging station. Each drone $u \in U$ is associated with a starting and an ending depot $D^S_u,D^E_u \in \VD$. Each arc $(i,j) \in A$ is associated with a travel time $t_{ij}$, a travel distance $d_{ij}$, a travel cost $c_{ij}$, e.g., required battery power, a handover number $h_{ij}$, and an expected outage duration $o_{ij}$. 

Let $T_u=\{T_u^1,T_u^2,\ldots,T_u^{n_C}\}$ be the feasible set of trips of drone $u \in U$, where each item in this set consists of arc sequences used in a particular trip, i.e., $T_u^k = \{(v_1,v_2),\ldots,(v_{{n_k}-1},v_{n_k}): v_1 = D^s_u, v_{n_k} = D^E_u,$ $(v_{i-1},v_i) \in A, i=2,\ldots,n_k\}$ denotes the $k^{\textrm{th}}$ trip of drone $u$ with $n_k$ arcs flown. Note that the maximum number of trips could be equal to the number of customers for a particular drone since each trip can only serve to a single customer. In case a drone has less trip than the number of customers, then, the corresponding trip is empty.

For a given trip $T_u^k$, we introduce the following equations:
\begin{flalign}
    H \left( T_u^k \right) & = \sum_{(i,j)\in T_u^k} h_{ij}, \label{eqn:H} \\
    O \left( T_u^k \right) & = \sum_{(i,j)\in T_u^k} o_{ij}, \label{eqn:O} \\
    M \left( T_u^k \right) & = \sum_{(i,j)\in T_u^k} d_{ij}, \label{eqn:M}
\end{flalign}

\noindent where $H$, $O$, and $M$ define the number of handover operations, the outage duration, and the total flight distance, respectively. Let $\Omega$ denote the set of all feasible trip sets for all drones. Then, given the handover and outage QoS thresholds $\Hmax$ and $\Omax$, respectively, the C-DDP is defined as:

\begin{flalign}
    \nonumber \textrm{minimize} \ & \quad \sum_{u \in U} \sum_{k \in \{1,\ldots,n_C \}} M \left( T_u^k \right) \\
    \textrm{subject to} & \quad H \left( T_u^k \right) \leq \Hmax, u \in U, k \in \{1,\ldots,n_C \}, \label{eqn:compact_const_1} \\
    & \quad O \left( T_u^k \right) \leq \Omax, u \in U, k \in \{1,\ldots,n_C \},   \label{eqn:compact_const_2} \\
    & \quad \mathbf{T} \in \Omega,   \label{eqn:compact_const_3}
\end{flalign}

\noindent where $\mathbf{T}$ is the set of all trips of all drones, i.e., $\mathbf{T}= \{T_{U_1},T_{U_2},\ldots,T_{U_{n_U}}\}$. The objective of the C-DDP is to minimize the total flight distance of trips while constraints \eqref{eqn:compact_const_1} and \eqref{eqn:compact_const_2} ensure that the number of required handover operations and expected outage duration of any trip do not exceed the QoS thresholds, respectively.

\begin{figure*}[!htb]
    \centering
    \begin{tabular}{cc}
        \subfloat[$\Hmax = \infty$, $\Omax = \infty$.]{
            \includegraphics[width=.45\textwidth]{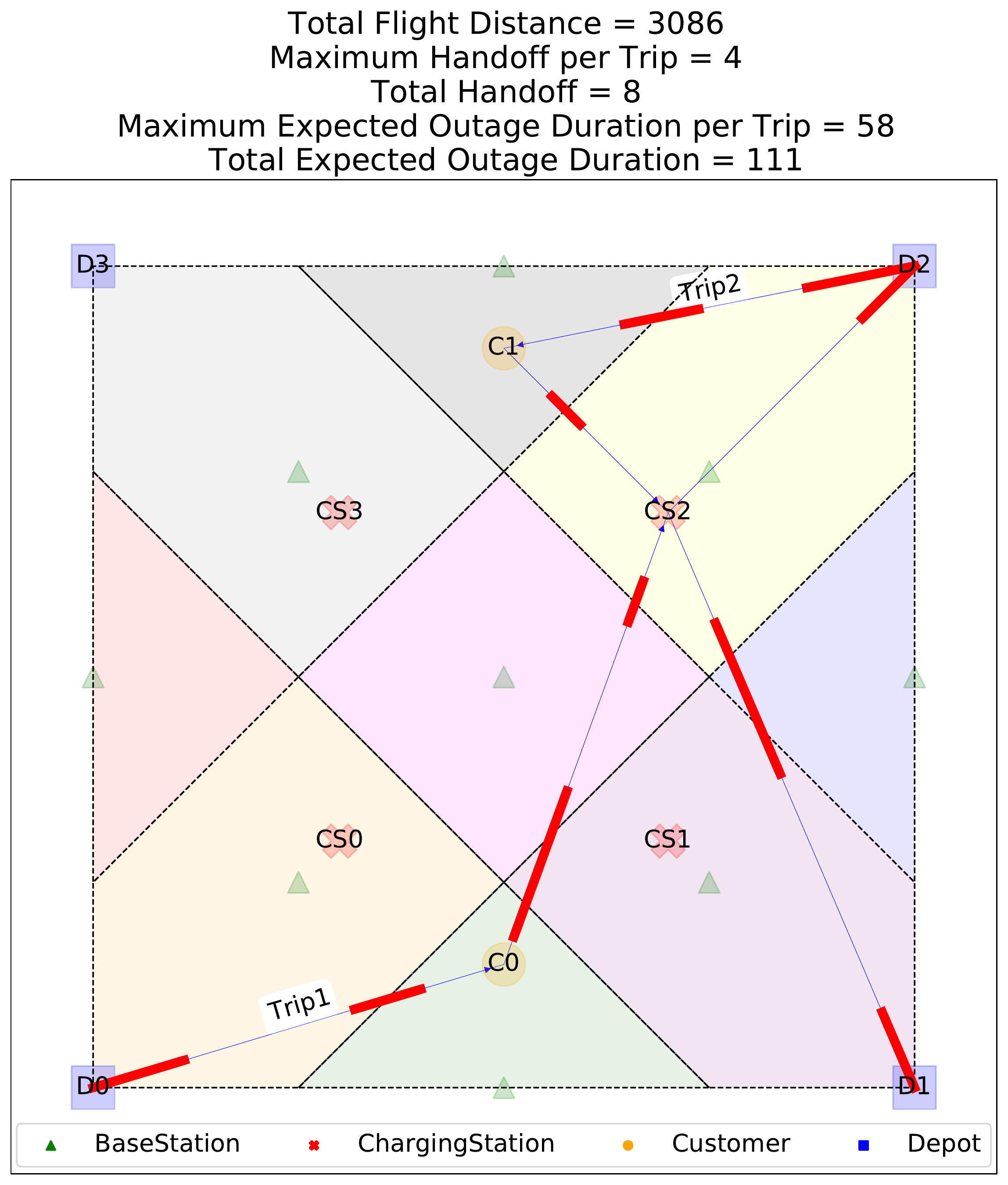}
            \label{fig:illustrative_example_default}
        }
        &
        \subfloat[$\Hmax = 3$, $\Omax = 45$.]{
            \includegraphics[width=.45\textwidth]{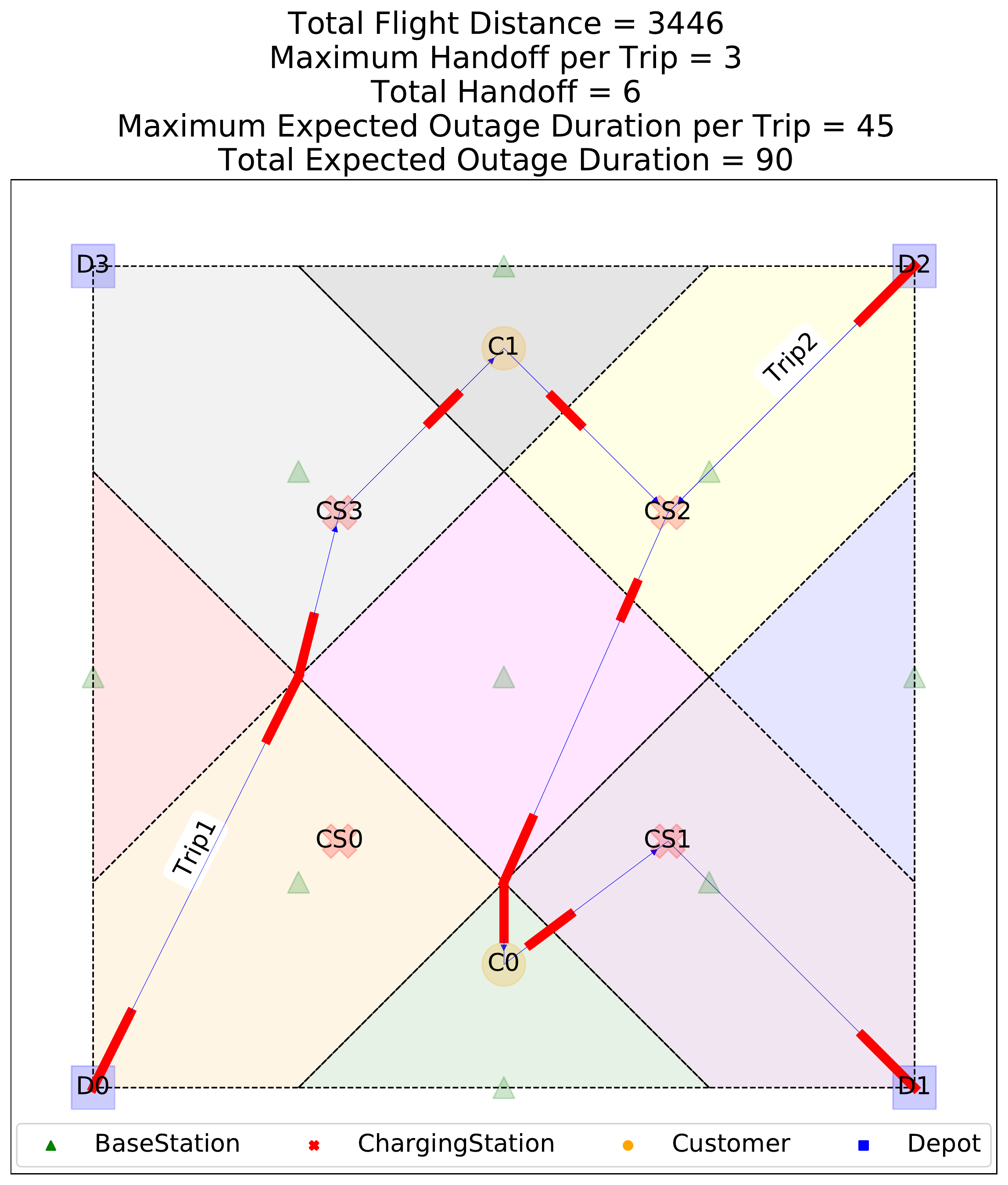}
            \label{fig:illustrative_example_3_10}
        }
        \\
        \subfloat[$\Hmax = 3$, $\Omax = 51$.]{
            \includegraphics[width=.45\textwidth]{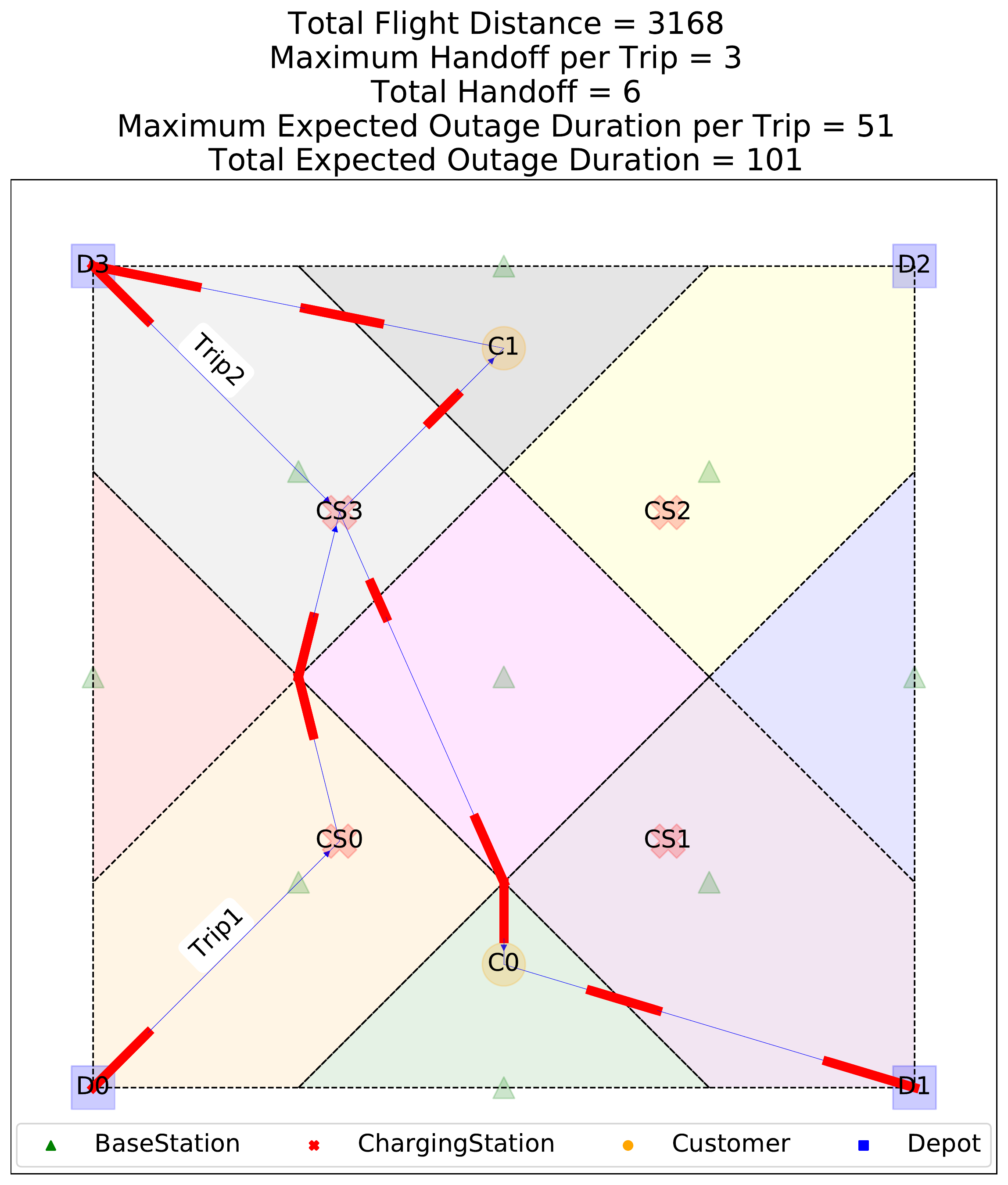}
            \label{fig:illustrative_example_3_20}
        }
        &
        \subfloat[$\Hmax = 3$, $\Omax = 56$.]{
            \includegraphics[width=.45\textwidth]{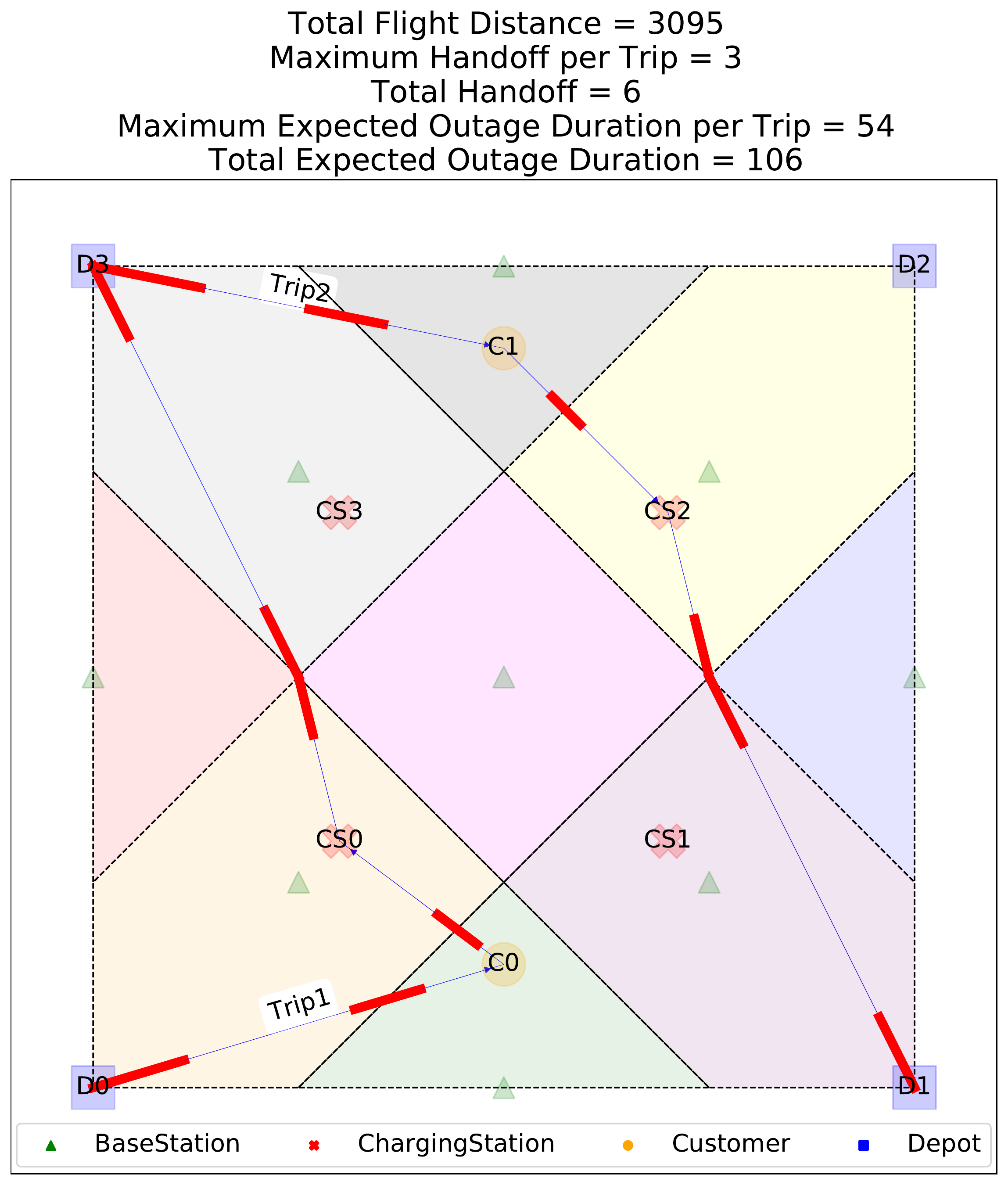}
            \label{fig:illustrative_example_3_30}
        }
    \end{tabular}
    \caption{\it \textbf{Illustrative example given $\mathbf{\Bar{\gamma}=2}$.} The thicker red parts of the trips indicate the expected outage. The top-left panel shows the solution for the default problem. The top-right, bottom-left, and bottom-right panels show the case in which the outage duration constraint is relaxed by 10\%, 20\%, and 30\%, respectively, while the handover constraint is set to its optimum.}
    \label{fig:illustrative_example}
\end{figure*}

\subsection{Illustrative Example}\label{subsec:illustrative_example}
In Fig.~\ref{fig:illustrative_example}, we illustrate the C-DDP on a $1000\times 1000$ square region with two customers (C) located at $[500,150]$ and $[500,900]$, four depots (D) located at $[0,0]$, $[1000,0]$, $[1000,1000]$, and $[0,1000]$, and four charging stations (CS) located at $[300,300]$, $[700,300]$, $[700,700]$, and $[300,700]$. The communication network includes 9 base stations located on equidistant diagonals and $\Bar{\gamma}$ is set to $2$. 

We assume that a single drone will deliver the customer packages starting from depot ``D0'' and ending at depot ``D1''. For simplicity, we ignore the time windows in this example. The different colored regions show the coverage areas in the communication network, and the dashed lines depict the boundaries of these areas. The top left panel shows the optimal solution for the default case, where the objective is to minimize the total flight distance and constraints~\eqref{eqn:compact_const_1}, \eqref{eqn:compact_const_2} are ignored. We have also solved the same problem by considering handover and outage objectives to determine the optimal communication performance. In particular, we set the objective function of C-DDP to minimize the maximum number of handovers, i.e., $\min_{u \in U}$ $\max_{k \in \{1,\ldots,n_C \}} H(T_u^k)$, and the maximum expected outage duration per trip, i.e., $\min_{u \in U} \max_{k \in \{1,\ldots,n_C \}} O(T_u^k)$, with no additional constraints. We obtained the optimal handover count and outage duration as 3 and 43, respectively.

We illustrate three different cases in Fig.~\ref{fig:illustrative_example}, where we set the objective function to minimize the total flight distance and consider the optimal handover and outage duration as hard constraints. Note that the problem becomes infeasible once we enforce both the optimal handover and outage constraints, i.e., $\Hmax = 3$ and $\Omax = 43$. Therefore, we illustrate the cases in which we relax the maximum outage duration by 10\% (top right), 20\% (bottom left), and 30\% (bottom right).

We can easily observe that integrating the communication constraints increases the total flight distance. In particular, the flight distance increases by 12.3\%, 11.7\%, and 6.7\% when at most 10\%, 20\%, and 30\% relaxation are allowed for the maximum expected outage duration, respectively. However, discarding the communication constraints would end up in a 34.8\% increase in the expected outage duration in the default case compared to the optimal expected outage duration per trip. Note that relaxing the handover constraint would not affect these solutions since only the outage constraint has become active in this particular example.

\subsection{MIP Formulation of the C-DDP}\label{subsec:mip}
Let $K=\{0,1,2,\ldots,n_C-1\}$ denote the set of trips for each drone. We introduce the decision variables $x_{ijuk} \in \{0,1\}$ to denote if arc $(i,j) \in A$ is flown by drone $u \in U$ on its $k^{\textrm{th}}$ trip, 0 otherwise; $p_{uk} \in \{0,1\}$ to denote if drone $u \in U$ operates trip $k \in K$, 0 otherwise; $0 \leq y_{iuk} \leq 1$ to denote the battery level of drone $u \in U$ at the time of arrival at node $i \in \VF$ on its $k^{\textrm{th}}$ trip; $s^L_{iuk} \geq 0$ and $s^A_{iuk} \geq 0$, respectively, to denote leaving and arrival time of drone $u \in U$ at depot $i \in \VD$ on its $k^{\textrm{th}}$ trip; and $s^V_{iuk} \geq 0$ to denote the arrival time of drone $u \in U$ at node $i\in \VF$ on its $k^{\textrm{th}}$ trip. Then, the MIP formulation of C-DDP can be defined as
{\allowdisplaybreaks
\begin{flalign}
    \nonumber \textrm{minimize} \ & \sum_{(i,j) \in A} \sum_{u \in U} \sum_{k \in K} x_{ijuk} d_{ij} & \\
    \textrm{subject to} \ & \sum_{(i,j) \in A} h_{ij} x_{ijuk} \leq \Hmax, u \in U, k \in K, \label{mip::cons_handover} \\
    & \sum_{(i,j) \in A} o_{ij} x_{ijuk} \leq \Omax, u \in U, k \in K, \label{mip::cons_outage} \\
    & \sum_{\substack{i \in \VF, \\ (i,j) \in A}} \sum_{u \in U} \sum_{k \in K} x_{ijuk} = 1, j \in \VC, \label{mip::cons_customer_visit} \\
    & \sum_{\substack{j \in \VF \\ (i,j) \in A}} x_{ijuk} - \sum_{\substack{j \in \VF, \\ (j,i) \in A}} x_{jiuk} = 0, i \in \VF \setminus \VD, u \in U, k \in K, \label{mip::cons_flow_balance} \\
    & p_{u0} \leq \sum_{\substack{j \in \VF \\ (D^S_u,j) \in A}} x_{D^S_u,j,u,0} , u \in U, \label{mip::cons_start_depot} \\
    & p_{uk} - p_{u,k+1} \leq \sum_{\substack{j \in \VF \\ (j,D^E_u) \in A}} x_{j,D^E_u,u,k} , u \in U, k \in K \setminus \{n_C-1\}, \label{mip::cons_end_depot_1} \\
    & p_{u,n-1} \leq \sum_{\substack{j \in \VF \\ (j,D^E_u) \in A}} x_{j,D^E_u,u,n-1} , u \in U, \label{mip::cons_end_depot_2} \\
    & x_{ijuk} \leq p_{uk} , (i,j) \in A, u \in U, k \in K, \label{mip::cons_arc_assignment} \\
    & p_{u,k+1} \leq p_{uk} , u \in U, k \in K \setminus \{n_C-1\}, \label{mip::cons_trip_sequence} \\
    & \sum_{i \in \VD} \sum_{\substack{j \in \VF \\ (i,j) \in A}} x_{ijuk} \leq 1, u \in U, k \in K, \label{mip::cons_only_one_leaving_depot} \\
    &  \sum_{i \in \VD} \sum_{\substack{j \in \VF \\ (j,i) \in A}} x_{jiuk} \leq 1, u \in U, k \in K, \label{mip::cons_only_one_arrival_depot} \\
    &  \sum_{i \in \VC} \sum_{\substack{j \in \VF \\ (j,i) \in A}} x_{jiuk} \leq 1, u \in U, k \in K, \label{mip::cons_only_one_customer} \\
    & \sum_{\substack{j \in \VF \\ (i,j) \in A}} x_{i,j,u,k+1} - \sum_{\substack{j \in \VF \\ (j,i) \in A}} x_{jiuk} \leq 0 , i \in \VD, u \in U, k \in K \setminus \{n_C-1\}, \label{mip::cons_depot_in_out} \\
    & x_{ijuk} - 1 \leq y_{juk} - (1 - c_{ij}) \leq 1 - x_{ijuk}, i \in \VCS, j \in \VF, (i,j) \in A, u \in U, k \in K, \label{mip::cons_battery_1} \\    
    & x_{ijuk} - 1 \leq y_{juk} - (y_{iuk} - c_{ij}) \leq 1 - x_{ijuk}, i \in \VF \setminus \VCS, j \in \VF, (i,j) \in A, u \in U, k \in K, \label{mip::cons_battery_2} \\
    & a_i - M \left( 1 - \sum_{\substack{j \in \VF \\ (j,i) \in A}} x_{jiuk} \right) \leq s^V_{iuk} \leq b_i + M \left( 1 - \sum_{\substack{j \in \VF \\ (j,i) \in A}} x_{jiuk} \right), i \in \VC, u \in U, k \in K, \label{mip::cons_time_window} \\
    & s^A_{iuk} - M(1 - p_{u,k+1}) \leq s^L_{i,u,k+1}, i \in \VD, u \in U, k \in K \setminus \{n_C-1\}, \label{mip::cons_trip_start_time} \\
    & s^L_{iuk} + w_i + t_{ij} - M(1 - x_{ijuk}) \leq s^V_{juk}, i \in \VD, j \in \VF \setminus \VD, (i,j) \in A, u \in U, k \in K, \label{mip::cons_first_node_visit_time} \\
    & s^V_{iuk} + w_i + t_{ij} - M(1 - x_{ijuk}) \leq s^V_{juk}, i,j \in \VF \setminus \VD, (i,j) \in A, u \in U, k \in K, \label{mip::cons_interim_node_visit_time} \\
    & s^V_{iuk} + w_i + t_{ij} - M(1 - x_{ijuk}) \leq s^A_{juk}, i \in \VF \setminus \VD, j \in \VD, (i,j) \in A, u \in U, k \in K, \label{mip::cons_trip_end_time_from_nondepot} \\
    & s^L_{iuk} + t_{ij} - M(1 - x_{ijuk}) \leq s^A_{juk}, i,j \in \VD, (i,j) \in A, u \in U, k \in K, \label{mip::cons_trip_end_time_from_depot} \\
    & x_{ijuk} \in \{0,1\}, (i,j) \in A, u \in U, k \in K, \label{mip::cons_domain_1} \\
    & p_{uk} \in \{0,1\}, u \in U, k \in K, \label{mip::cons_domain_2} \\
    & 0 \leq y_{iuk} \leq 1, i \in \VF, u \in U, k \in K, \label{mip::cons_domain_3} \\
    & 0 \leq s^L_{iuk}, s^A_{iuk}, i \in \VD, u \in U, k \in K, \label{mip::cons_domain_4} \\
    & 0 \leq s^V_{iuk}, i \in \VF \setminus \VD, u \in U, k \in K. \label{mip::cons_domain_5}
\end{flalign}
}

The objective function minimizes the total flight distance. Constraints~\eqref{mip::cons_handover} and \eqref{mip::cons_outage} ensure that the entire delivery operation is fulfilled without exceeding the handover and outage thresholds, respectively. Constraints~\eqref{mip::cons_customer_visit} state that each customer is visited once by one drone. Constraints~\eqref{mip::cons_flow_balance} define the flow balance at each vertex. Constraints~\eqref{mip::cons_start_depot}-\eqref{mip::cons_end_depot_2} ensure that a drone starts the workday from its starting depot and ends the day at its ending depot. In particular, \eqref{mip::cons_start_depot} enforces to have an outgoing arc from starting depot of a drone if that drone is assigned its first trip, while \eqref{mip::cons_end_depot_1} and \eqref{mip::cons_end_depot_2} enforce to have an incoming arc to the ending depot of a drone if that drone is not assigned a new trip after operating its last trip. Constraints~\eqref{mip::cons_arc_assignment} ensure that arcs can be flown by a particular drone on a particular trip if that trip is assigned to that drone. Constraints~\eqref{mip::cons_trip_sequence} organize the trip sequences. Constraints~\eqref{mip::cons_only_one_leaving_depot}-\eqref{mip::cons_only_one_customer} require that each trip can have at most one leaving depot, one arrival depot, and a customer. Constraints~\eqref{mip::cons_depot_in_out} ensure that a drone can start its next trip from a particular depot if and only if the previous trip ends at that particular depot. Constraints~\eqref{mip::cons_battery_1}-\eqref{mip::cons_battery_2} determine the battery levels.

Constraints~\eqref{mip::cons_time_window}-\eqref{mip::cons_trip_end_time_from_depot} state the time relationships. In particular, Constraints~\eqref{mip::cons_time_window} ensure that a delivery is fulfilled within the time window of a particular customer. Constraints~\eqref{mip::cons_trip_start_time} ensure that the leaving time of the next trip for a drone cannot be earlier than the arrival time of the previous trip. Constraints~\eqref{mip::cons_first_node_visit_time} set the visiting time of the first non-depot node in a trip to the leaving time from the starting depot of the trip plus the operation time at the depot plus the travel time from depot to the first node. Constraints~\eqref{mip::cons_interim_node_visit_time} determine the visiting time of a non-depot node as the visiting time of the previous non-depot node plus operation time spent in the previous node plus the travel time from the previous node. The operation times in these constraints correspond to loading, unloading, and charging times for depots, customers, and charging nodes, respectively. Constraints~\eqref{mip::cons_trip_end_time_from_nondepot} set the arrival time to a depot node from a non-depot node to the visiting time of the last non-depot node plus the operation time at that last node plus the travel time from that last node to the depot node. Constraints~\eqref{mip::cons_trip_end_time_from_depot} set the arrival time to a depot node from another depot node to the leaving time of the outgoing depot plus the travel time between depots. Note that there is no operation time in this constraint since we assume that flying between two depot nodes does not require any loading/unloading operation. Finally, constraints~\eqref{mip::cons_domain_1}-\eqref{mip::cons_domain_5} define the domains of the variables.

\section{Solution Approach}\label{sec:solution_approach}

The MIP formulation moves out of our computational reach for large problem instances. Therefore, we propose a \emph{Genetic Algorithm} (GA) to effectively solve larger instances. We have developed an individual representation for the GA and implemented it via the \emph{pymoo} package \cite{blank2020}, which incorporates different evolutionary algorithms for single- and multi-objective problems.

An individual in GA defines how a solution, i.e., a set of trips of our drone fleet, can be stated. We can define a vector of nodes over which a drone performs its trips in our context. Then, we aggregate all vectors to define an individual of the GA. Note that the size of a vector is not known a priori. However, we fix this size to the number of customers times the number of flyable nodes. This is because each customer must be visited once, and a customer is guaranteed to be visited by a drone after all flyable nodes in the network are traversed in the worst case.

\begin{figure*}[!htb]
    \centering
    \begin{tabular}{c}
        \subfloat[Label encoding.]{
            \includegraphics[width=.7\textwidth]{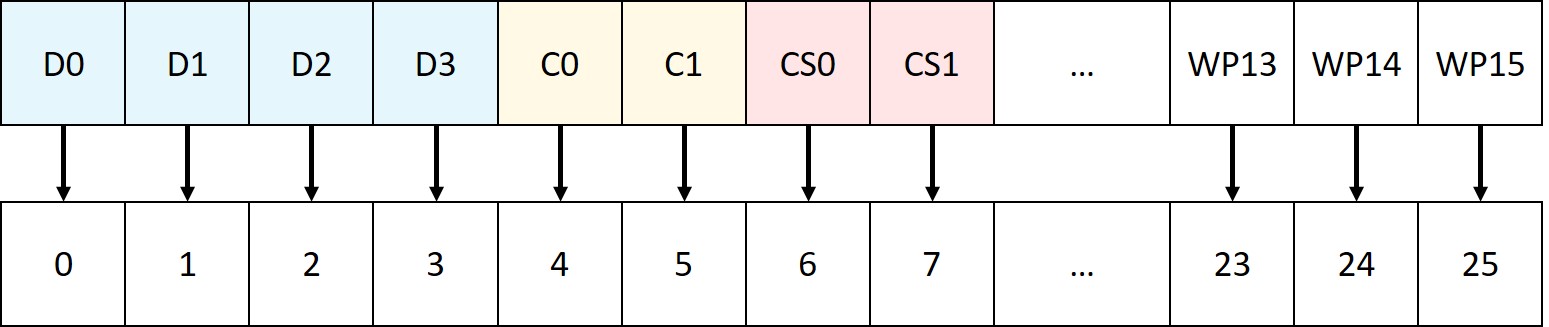}
            \label{fig:individual_representation_a}
        }
        \\
        \subfloat[Individual representations.]{
            \includegraphics[width=.9\linewidth]{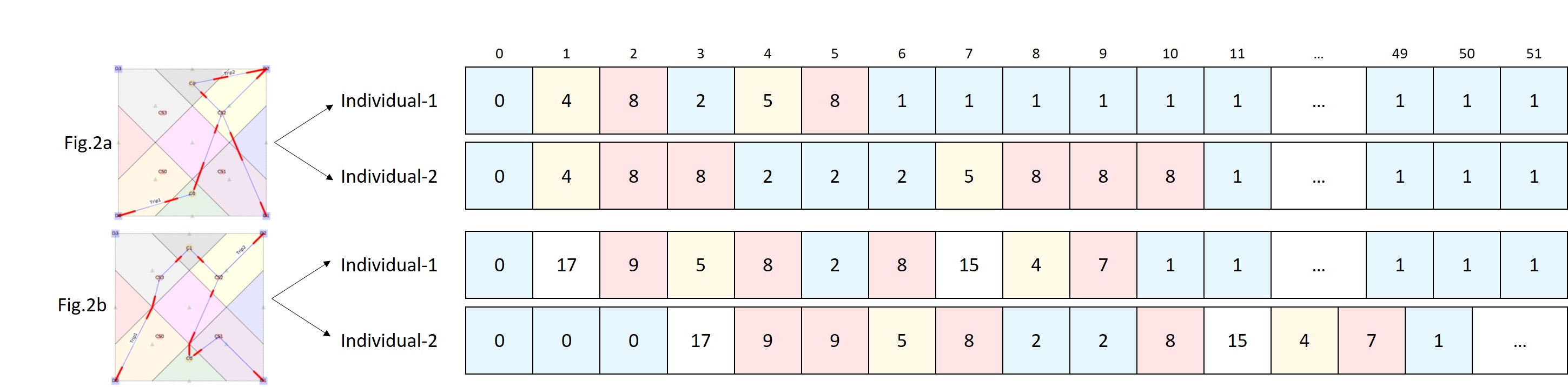}
            \label{fig:individual_representation_b}
        }
    \end{tabular}
    \caption{\it \textbf{Individual representation.}The top panel shows the label encoding of flyable nodes including depots, customers, charging stations, and waypoints. The bottom panel shows the individual representations (recovered at the top and unrecovered at the bottom) of trips for the top two panels in the illustrative example in Fig.~\ref{fig:illustrative_example}.}
    \label{fig:individual_representation}
\end{figure*}

Fig.~\ref{fig:individual_representation} demonstrates how we represent the individuals. In particular, we first encode every flyable node in the network using labels. Fig.~\ref{fig:individual_representation_a} shows the encoding for the illustrative example presented in Fig.~\ref{fig:illustrative_example}. Recall that we have 4 depots, 4 charging stations, and 2 customers in this network. Moreover, we have 9 base stations. To incorporate possible \emph{Waypoints} (WPs), we include the Voronoi points of base stations, i.e., the intersection points of coverage areas, as flyable nodes. We also add the intersections with boundaries of the service area and coverage areas as WPs. As a result, we have 26 labels for this particular example.

In Fig.~\ref{fig:individual_representation_b}, two individuals are shown for Fig.~\ref{fig:illustrative_example_default} and \ref{fig:illustrative_example_3_10}, respectively. Since we have 26 flyable nodes in the network after adding the WPs, each individual has 52 components for the single drone case. This number would increase to 104 if we had two drones, where the first 52 components show the trips of the first drone, and the remaining 52 components show the trips for the second drone. Note that the same label may be sequentially seen in a particular individual. In such cases, we recover the routes by removing repetitions and considering only one label at a time. This recovery operation can be seen in Fig.~\ref{fig:individual_representation_b}. The individual at the top (Individual-1) shows the recovered representation of the individual at the bottom (Individual-2) for the same solution.

Based on the individual representation we have introduced, we can use the standard genetic operators like one-point or random crossover, polynomial mutation, and random sampling, which allows for a straightforward implementation of the algorithm. Moreover, the constraints can be handled in \emph{pymoo} such that violations can be penalized. We use the \emph{ConstraintsAsPenalty} method with its default parameters to handle constraints.

\section{Computational Study}\label{sec:computational_study}
In this section, we present the results of our computational experiments to validate and assess our proposed model as well as the performance of the GA. We conducted all experiments on a desktop computer with an Intel i-7@3.50 GHz processor, 16 GB RAM, and Windows 10 operating system. The MIP models were solved with Gurobi 9.5 \cite{gurobi}.

\subsection{Test Bed}\label{subsec:results_synthetic_data_generation}
Since there are no publicly available benchmark instances for the C-DDP, we have generated test instances to conduct numerical experiments\footnote{Instances are available at \url{https://github.com/cihantugrulcicek/CDDP}.}. We have built our instances using a $5 \times 5~\mathrm{km}$ zone for a working day of 8 hours (28,800 seconds). We have located the depot and CSs for all instances such that there exists one depot every $2~\mathrm{km}$ and one CS every $1~\mathrm{km}$. Then, we have generated the communication network with respect to two different spatial distributions. In particular, we have first created a hexagonal grid with $1~\mathrm{km}$ radius and located one GBS at the center of each hexagon for the uniform case. We have generated another network setting where we have perturbed the center location of hexagons such that the coverage areas overlap asymmetrically. 

To determine the customer locations, we use two different distributions. For the first case, the locations are determined according to \emph{Uniform} (U) distribution. For the second case, we use the Poisson Point process, where we first generate multiple hotpoints, and then determine customer locations around these hotpoints. The number of customers around each hotpoint is also determined based on a discrete uniform distribution. 

The time window for visiting a customer is determined as follows: we first determine the earliest time that a customer can be visited according to its travel time to the closest and the farthest depot. More precisely, let $D^\prime$ and $D^{\prime\prime}$ be the closest and farthest depots to customer $i \in \VC$, respectively. We set the center of time window by generating a uniform random number from $\Bar{a}_i = \mathrm{U}(t_{D^\prime ,i}, 8 - t_{i,D^{\prime\prime}})$. We then generate another uniform random number as the width of the time window and modify the exact time window with respect to this width. Note that in case of having an earliest time less than $t_{D^\prime ,i}$ or a latest time greater than ($8 - t_{i,D^{\prime\prime}}$), we adjust the time window such that the earliest time is set to $t_{D^\prime ,i}$ and the latest time is set to ($8 - t_{i,D^{\prime\prime}}$). To test the impact of time windows, we generate the width from two different Uniform distributions, i.e., we draw a random integer number from $\mathrm{U}(2,8)$ for the loose case, and $\mathrm{U}(1,4)$ for the tight case. 

In all instances, the number of drones is determined based on one drone per 25 customers. The starting and ending depots of drones are randomly determined, assuring that each depot hosts at least one drone at the beginning and end of the planning horizon if the number of drones is greater than the number of depots. The communication parameters are set to dense-urban settings provided by \cite{boryaliniz2016}: $\alpha_1=12.08$, $\alpha_2=0.11$, $\alpha_3=2.5$, $\mu_\mathrm{LoS}=1.6$, $\mu_\mathrm{NLoS}=23$, $\sigma^2=-173~\mathrm{dB}$, $P_i=46~\mathrm{dBm}$ for all $i\in\VCN$.

To succinctly describe the instances, we use a three-digit notation. ``U'' or ``P'' in the first digit denotes Uniform and Perturbed communication network, ``U'' or ``P'' in the second digit denotes Uniform and Poisson Process customer distribution,  and ``L'' or ``T'' in the last digit denotes loose and tight time windows, e.g., ``PUT'' indicates the setting with the Perturbed communication network, Uniform delivery network and Loose time windows.

\subsection{Results}

We set the Gurobi parameters to their defaults except that the time limit has been set to 3600 seconds. In GA, we use standard integer binary crossover (\emph{int\_sbx}), polynomial mutation (\emph{int\_pm}), random sampling (\emph{int\_random}) operators, and all other parameters such as population size and number of generations was set to \emph{pymoo} defaults. The GA is terminated whenever the run time exceeds 3600 seconds, or the number of generations exceeds ten times the number of customers ($10n_C$) in a particular instance, e.g., 500 generations for an instance with 50 customers.

We first compare the GA's solution performance (optimality gap) against Gurobi. In particular, let $M^\mathrm{best}$ denote the best feasible objective value, $M^\mathrm{G}$ denote the best bound reported by Gurobi, and $M^\mathrm{GA}$ denote the best feasible GA objective value. Then, $(M^\mathrm{best} - M^\mathrm{G}) \slash M^\mathrm{G}$ and $(M^\mathrm{GA} - M^\mathrm{G})  \slash M^\mathrm{G}$ are the optimality gaps of Gurobi and GA, respectively.

\begin{table}[!b]
    \small
    \centering
    \caption{Optimality gaps of Gurobi and the GA.}
    \begin{tabular}{llrrrrrrcrrrr}
    \toprule
    && \multicolumn{5}{c}{Gurobi} & & \multicolumn{5}{c}{GA} \\ \cline{3-7} \cline{9-13}
    \multicolumn{1}{c}{Setting} & \multicolumn{1}{c}{$n_C$} &  \multicolumn{1}{c}{Default}  & \multicolumn{1}{c}{(20,20)} & \multicolumn{1}{c}{(20,10)} & \multicolumn{1}{c}{(10,20)} & \multicolumn{1}{c}{(10,10)} & & \multicolumn{1}{c}{Default} & \multicolumn{1}{c}{(20,20)} & \multicolumn{1}{c}{(20,10)} & \multicolumn{1}{c}{(10,20)} & \multicolumn{1}{c}{(10,10)}  \\
    \midrule
    PPL & 5  &   0.019 &   0.026 &   0.039 &   0.036 &   0.045 & &  0.018 &   0.026 &   0.038 &   0.036 &   0.043 \\
    & 10 &   0.020 &   0.027 &   0.040 &   0.039 &   0.043 &  & 0.020 &   0.027 &   0.041 &   0.039 &   0.044 \\
    & 20 &   0.021 &   0.026 &   0.043 &   0.041 &   0.044 &  & 0.020 &   0.027 &   0.045 &   0.041 &   0.045 \\
    \midrule
    PPT & 5  &   0.018 &   0.000 &   0.000 &   0.000 &   0.000 &  & 0.018 &   0.000 &   0.000 &   0.000 &   0.000 \\
    & 10 &   0.019 &   0.000 &   0.000 &   0.025 &   0.000 & &  0.019 &   0.000 &   0.000 &   0.026 &   0.000 \\
    & 20 &   0.019 &   0.042 &   0.000 &   0.065 &   0.005 &  & 0.019 &   0.042 &   0.000 &   0.067 &   0.021 \\
    \midrule
    PUL & 5  &   0.000 &   0.000 &   0.000 &   0.000 &   0.000 &  & 0.000 &   0.000 &   0.000 &   0.000 &   0.000 \\
    & 10 &   0.015 &   0.000 &   0.000 &   0.000 &   0.000 &  & 0.014 &   0.000 &   0.000 &   0.000 &   0.000 \\
    & 20 &   0.000 &   0.000 &   0.000 &   0.000 &   0.125 &  & 0.000 &   0.000 &   0.000 &   0.000 &   0.121 \\
    \midrule
    PUT & 5  &   0.154 &   0.080 &   0.046 &   0.080 &   0.046 &  & 0.154 &   0.080 &   0.045 &   0.082 &   0.044 \\
    & 10 &   0.165 &   0.085 &   0.050 &   0.081 &   0.050 &  & 0.160 &   0.087 &   0.050 &   0.079 &   0.048 \\
    & 20 &   0.177 &   0.088 &   0.048 &   0.087 &   0.051 &  & 0.173 &   0.084 &   0.050 &   0.088 &   0.051 \\
    \midrule
    UPL & 5  &   0.087 &   0.080 &   0.087 &   0.081 &   0.087 &  & 0.085 &   0.080 &   0.089 &   0.077 &   0.088 \\
    & 10 &   0.094 &   0.087 &   0.083 &   0.077 &   0.085 &  & 0.095 &   0.084 &   0.085 &   0.074 &   0.081 \\
    & 20 &   0.093 &   0.086 &   0.079 &   0.074 &   0.087 &  & 0.092 &   0.089 &   0.082 &   0.071 &   0.085 \\
    \midrule
    UPT & 5  &   0.144 &   0.154 &   0.175 &   0.154 &   0.175 &  & 0.141 &   0.156 &   0.170 &   0.146 &   0.166 \\
    & 10 &   0.145 &   0.153 &   0.175 &   0.153 &   0.166 & &  0.148 &   0.156 &   0.177 &   0.147 &   0.169 \\
    & 20 &   0.140 &   0.166 &   0.167 &   0.145 &   0.173 &  & 0.133 &   0.160 &   0.162 &   0.146 &   0.169 \\
    \midrule
    UUL & 5  &   0.115 &   0.121 &   0.080 &   0.121 &   0.079 & &  0.115 &   0.116 &   0.080 &   0.117 &   0.080 \\
    & 10 &   0.110 &   0.124 &   0.076 &   0.117 &   0.087 &  & 0.105 &   0.120 &   0.077 &   0.112 &   0.086 \\
    & 20 &   0.107 &   0.133 &   0.082 &   0.113 &   0.088 &  & 0.107 &   0.134 &   0.082 &   0.111 &   0.086 \\
    \midrule
    UUT & 5  &   0.168 &   0.155 &   0.146 &   0.121 &   0.137 & &  0.170 &   0.149 &   0.139 &   0.119 &   0.139 \\
    & 10 &   0.162 &   0.152 &   0.151 &   0.130 &   0.146 & &  0.166 &   0.152 &   0.154 &   0.131 &   0.150 \\
    & 20 &   0.165 &   0.161 &   0.145 &   0.130 &   0.149 &  & 0.166 &   0.156 &   0.138 &   0.128 &   0.148 \\
    \bottomrule
    \end{tabular}
    \label{tab:gurobi_vs_ga}
\end{table}

\begin{figure}[!t]
    \centering
    \includegraphics[width=.8\linewidth]{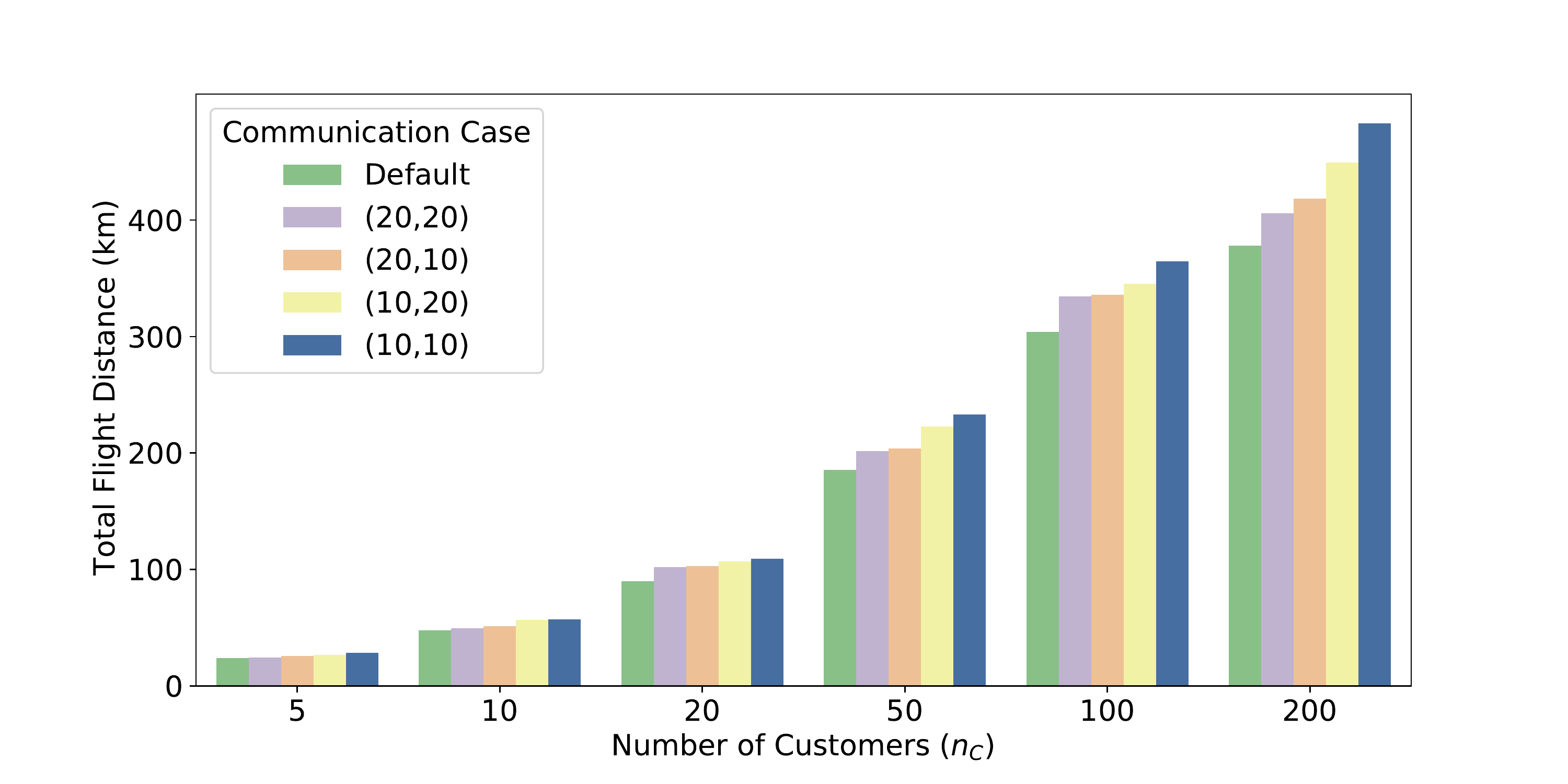}
    \caption{{\it \textbf{The change in total flight distance with respect to different levels of communication constraints.}} }
    \label{fig:communication_impact}
\end{figure}

\begin{figure}[!b]
    \centering
    \begin{tabular}{c}
        \subfloat[Maximum Number of Handovers.]{
            \includegraphics[width=.5\linewidth]{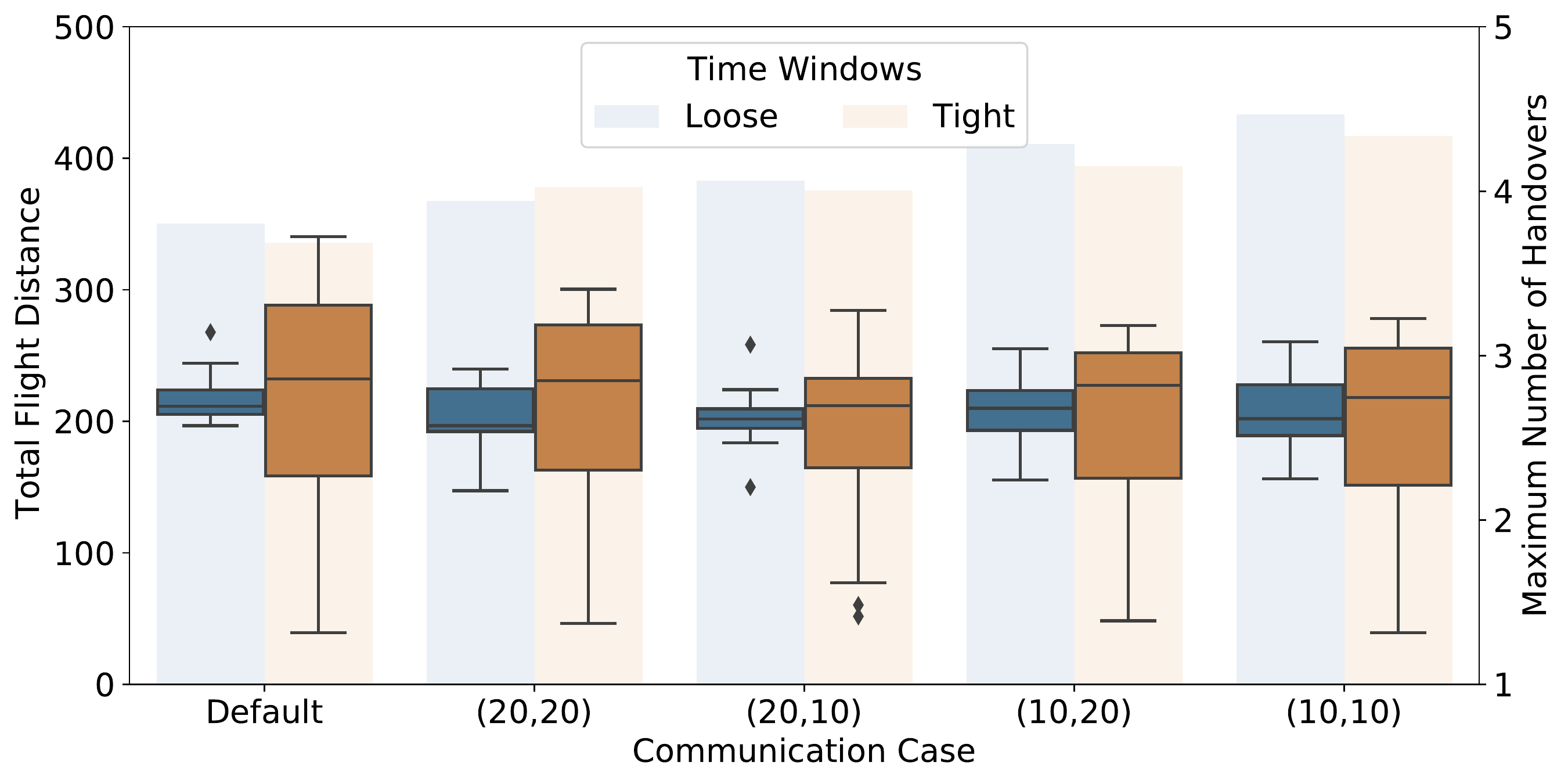}
            \label{fig:time_windows_handover}
        }
        \\
        \subfloat[Maximum Expected Outage.]{
            \includegraphics[width=.5\linewidth]{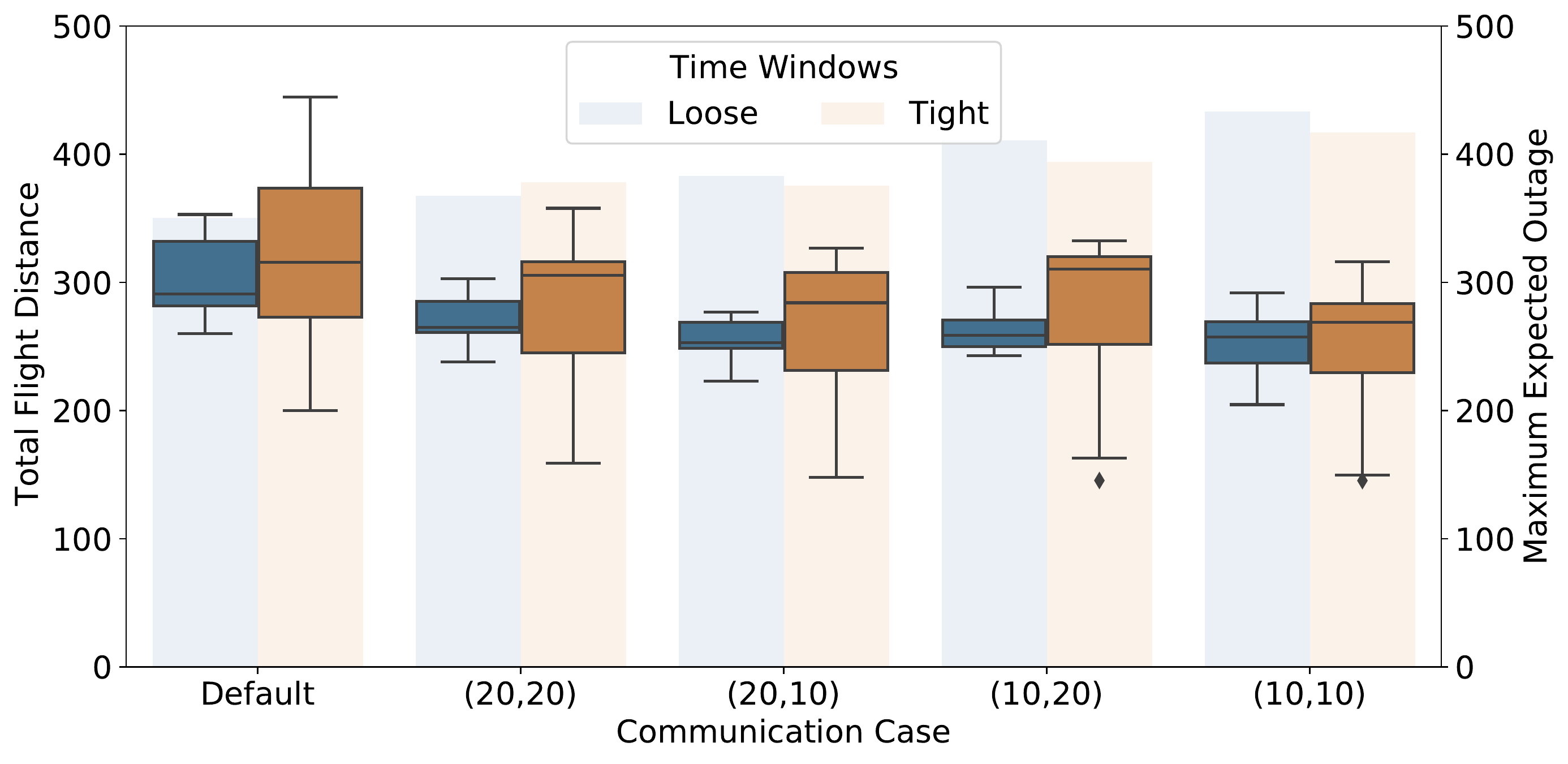}
            \label{fig:time_windows_outage}
        }
    \end{tabular}
    \caption{{\it \textbf{The impact of time windows on (a) maximum number of handovers and (b) maximum expected outage per trip.}} The bars and boxes indicate the total flight distance and the corresponding communication measure, respectively.}
    \label{fig:time_windows}
\end{figure}

\begin{figure}[!b]
    \centering
    \begin{tabular}{c}
        \subfloat[Maximum Number of Handovers.]{
            \includegraphics[width=.5\linewidth]{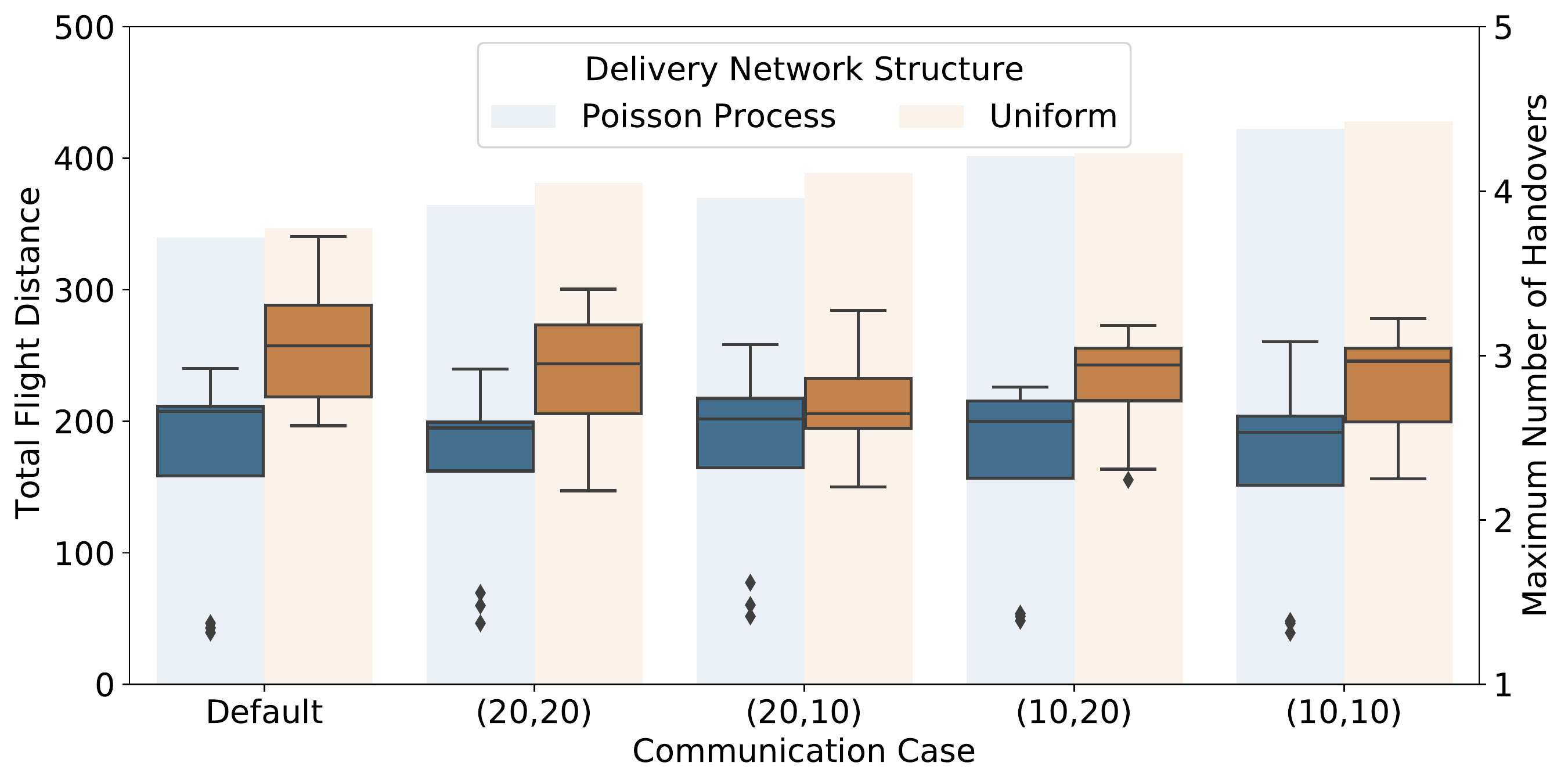}
            \label{fig:clustering_handover}
        }
        \\
        \subfloat[Maximum Expected Outage.]{
            \includegraphics[width=.5\linewidth]{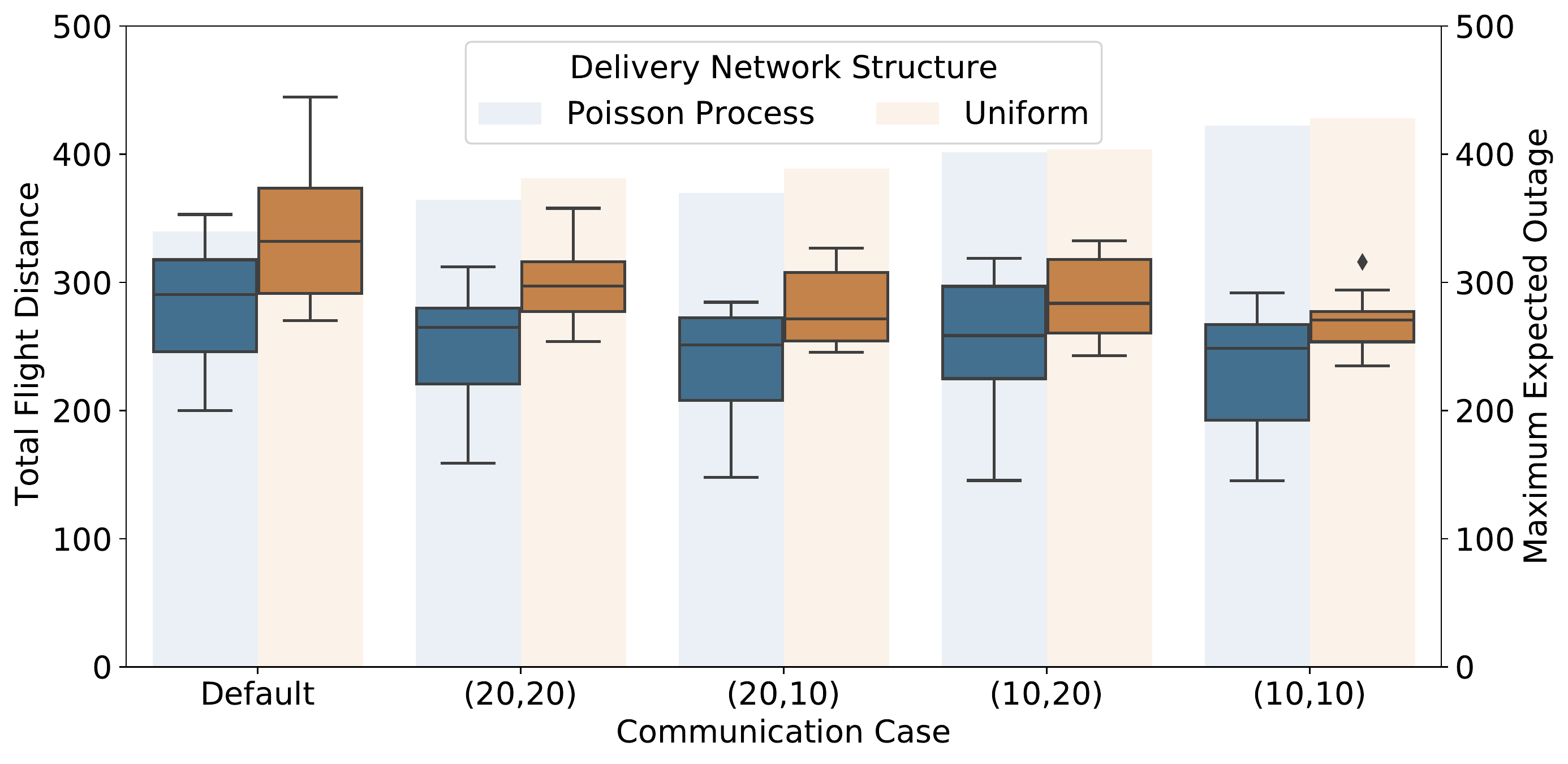}
            \label{fig:clustering_outage}
        }
    \end{tabular}
    \caption{{\it \textbf{The impact of customer locations on (a) maximum number of handovers and (b) maximum expected outage per trip.}} The bars and boxes indicate the total flight distance and the corresponding communication measure, respectively.}
    \label{fig:clustering}
\end{figure}

Table~\ref{tab:gurobi_vs_ga} presents the optimality gaps of Gurobi and the GA for different settings. Note that Gurobi could not find a feasible solution within the time limit for instances with more than 20 customers. We first solve each instance to find singular optimal values of maximum handover and expected outage. Then, we solve the same instances to find the optimal flight distances where we include the communication performance as hard constraints based on these singular optimal values. The columns in Table~\ref{tab:gurobi_vs_ga} denote how much we have relaxed each communication constraint. For example, ``(20,10)'' indicates that $\Hmax$ and $\Omax$ are set to 20\% and 10\% higher than their singular optimal values, respectively.

GA achieves approximately the same gap values with Gurobi. The worst and the best performances occur for the PPT setting with 0.11\% higher objective value and for the UPT setting with 0.26\% lower objective value, respectively. However, the solution times of the GA were more than 10 times shorter than Gurobi on the average, and more than 23 times shorter for instances with 20 customers. The average CPU time of Gurobi was 760 seconds while the GA has terminated under 70 seconds. This result confirms that the GA is a better option for solving larger instances.

Due to the limited computational reach of Gurobi in solving large instances, the remaining results are based solely on the replication averages of the GA. Fig.~\ref{fig:communication_impact} shows the impact of communication constraints on the total flight distance. The simulation results confirm our theoretical expectation that the lower the relaxation of a communication measure is, the higher the total flight distance. More precisely, the average total flight distance increases by 2.9\%, 3.6\%, 6.2\%, and 7.1\% for cases (20,20), (20,10), (10,20), and (10,10), respectively. On the other hand, the average number of handover activities decreases by 3.8\%, 6.4\%, 17.3\%, and 23.9\%, while the expected outage duration decreases by 11.3\%, 19.5\%, 26.5\%, and 31.2\%, respectively. This supports our claim that ignoring the communication constraints would result in operational disruption risk, which can be easily mitigated by sacrificing slightly from total flight distance, due the incorporation of $\Hmax$ and $\Omax$ into the formulation.

Fig.~\ref{fig:time_windows} and \ref{fig:clustering} illustrate the impact of time windows and clustering of customer requests on the total flight distance, number of handovers, and expected outage duration. The transparent bars in the background denote the average total flight distances, and the boxes indicate the corresponding communication measure. The type of time windows is not observed to have a significant impact on the total flight distance, the number of handovers, and expected outage. However, the delivery network structure significantly affects the communication performance and has a slight impact on the flight distance for all communication cases. In particular, the average number of handovers and expected outage duration decreases by 18\% and 16\% for instances with clustered customers. The clustered requests have also decreased the total flight distance by 3\% as the trips are likely to be organized from the closest depot to a hotpoint rather than flying to different depots in subsequent trips. This trip strategy requires fewer handover activities and consequently decreases the expected outage duration.

\section{Conclusions}\label{sec:conclusions}

In this study, we have addressed an emerging problem in last-mile delivery networks which we believe will attract more attention in the near future. In particular, we have integrated two significant communication performance measures, handover and expected outage duration, to the classical multi-depot multi-trip DDP. We have developed a new mathematical formulation, which can be solved by the off-the-shelf MIP solver Gurobi for instances with up to 20 customers. To solve larger problems, we have implemented a GA based on a novel individual representation that can be used in any kind of multi-depot multi-trip vehicle routing problem. Our numerical study has shown that the GA can outperform Gurobi in small-size instances and finds efficient solutions for larger sized instances with up to 200 customers.



\end{document}